\documentclass[a4paper,10pt]{amsart}

\usepackage[all]{xy}
\usepackage{amssymb}
\usepackage{amsmath}
\usepackage{amsthm}
\usepackage[english]{babel}

\newcommand{\ssection}{\vspace{0.2cm}\stepcounter{theorem}\noindent(\thetheorem) }

\makeatletter
\def\mylabel#1#2{\protected@edef\@currentlabel
       {\csname p@#1\endcsname\csname the#1\endcsname}
\label{#2}}
\makeatother

\newtheorem{theorem}{Theorem}[section]

\newtheorem{e-proposition}[theorem]{Proposition}

\newtheorem{e-definition}[theorem]{Definition\rm}
\newcommand\gl{{\mathbf{GL}}}
\renewcommand\O{\mathcal{O}}
\newcommand\git{/\!\!/}
\newcommand\tr{\mathrm{tr}}

\author{Olivier Serman}
\title{Local structure of $\mathbf{\mathcal{SU}_C(3)}$ for a curve of genus $\mathbf{2}$}
\address{Laboratoire J.-A.Dieudonn\'e \\ UMR 6621 du CNRS \\ Universit\'e de Nice \\ Parc Valrose \\ F-06108 Nice Cedex 02}
\email{serman@unice.fr}

\begin{document}
\maketitle

\begin{abstract}
% Text of abstract in English
The aim of this note is to give a precise description of the local structure of the moduli space $\mathcal{SU}_C(3)$ of rank $3$ vector bundles over a curve $C$ of genus $2$, which is in particular shown to be a local complete intersection. This allows us to investigate the local structure of the branch locus of the theta map, whose dual is known to be the Coble cubic in $\mathbb PH^0(J^1_C,3\Theta)$.
\end{abstract}

\section{Introduction}

Let $C$ be a smooth irreducible projective curve of genus $2$ over an algebraically closed field $k$ of characteristic zero, and let $\mathcal{SU}_C(3)$ be the moduli space of rank $3$ vector bundles over $C$ with trivial determinant. Laszlo began to investigate the local structure of this moduli space in \cite[V]{localstr}: Luna's \'etale slice theorem provides a way to compute the completed local ring at any point of $\mathcal{SU}_C(3)$ as GIT quotients of affine spaces, but, as soon as the isotropy group gets too bad, this leads to a quite intricate calculation. By translating this situation in terms of representations of quivers, we managed to work out the local structure at any point of $\mathcal{SU}_C(3)$. We have in particular obtained the following result:
\begin{theorem}
The moduli space of rank $3$ vector bundles over a curve of genus $2$ is a local complete intersection. 
\end{theorem}

As we have already seen in \cite{orth}, the notion of representations of quivers appears to be really helpful to understand the quotients given by Luna's result. Although it may not be clear in this note, where we could have given direct proofs avoiding such considerations, this quiver setting was the very basic point which led to generating sets for the coordinate rings of the quotients.

Let now $\Theta$ be the canonical Theta divisor on the variety $J^{1}$ which parametrizes line bundles of degree $1$ on $C$. It is known for long that the theta map $\theta \colon \mathcal{SU}_C(3) \to |3\Theta|$ is a double covering. Ortega has shown in \cite{angela} that its branch locus $\mathcal S \subset |3\Theta|$ is a sextic hypersurface which is the dual of the Coble cubic $\mathcal C \subset |3\Theta|^\ast$, where the Coble cubic is the unique cubic in $|3\Theta|^\ast$ which is singular along $J^1 \buildrel{|3\Theta|}\over{\longrightarrow}|3\Theta|^\ast$ (note that a different proof of this statement has been given by Nguy$\tilde{\hat{\text{e}}}$n in \cite{minh}). The last part of this paper is devoted to the local structure of the sextic $\mathcal S$.

\section{Local structure of $\mathcal{SU}_C(3)$} 

The starting point of the local study of moduli spaces of vector bundles, which follows from Luna's slice  theorem, can be found in \cite[II]{localstr}: it states that, at a closed point representing a polystable bundle $E$, the moduli space $\mathcal{SU}_C(r)$ of rank $r$ vector bundles with trivial determinant is \'etale locally isomorphic to the quotient $\mathrm{Ext}^1(E,E)_0 \git \mathrm{Aut}(E)$ at the origin, where $\mathrm{Ext}^1(E,E)_0$ denotes the kernel of $\tr \colon \mathrm{Ext}^1(E,E) \to H^1(C,\O_C)$.

We thus have to understand the ring of invariants of the polynomial algebra $k[\mathrm{Ext}^1(E,E)_0]=\mathrm{Sym}({\mathrm{Ext}^1(E,E)_0}^\ast)$ under the action of $\mathrm{Aut}(E)$. As a polystable bundle $E$ can be written 
\begin{align} \label{polystable}E=\bigoplus_{i=1}^s E_i \otimes V_i,\end{align}
\noindent where the $E_i$'s are mutually non-isomorphic stable bundles (of rank $r_i$ and degree $0$), and the $V_i$'s are vector spaces (of dimension $\rho_i$). Through this splitting our data become
\begin{align} \label{Ext}\mathrm{Ext}^1(E,E)=\bigoplus_{i,j} \mathrm{Ext}^1(E_i,E_j)\otimes \mathrm{Hom}(V_i,V_j),\end{align}
\noindent endowed with an operation of $\displaystyle{\mathrm{Aut}(E)=\prod_i \gl(V_i)}$ coming from the natural actions of $\gl(V_i) \times \gl(V_j)$ on $\mathrm{Hom}(V_i,V_j)$. 

We recognize here the setting of representations of quivers (see \cite{LBP}): consider indeed the quiver $Q$ with $s$ vertices $1,\ldots,s$, and $\dim\mathrm{Ext}^1(E_i,E_j)$ arrows from $i$ to $j$, and define $\alpha \in \mathbb N^s$ by $\alpha_i=\rho_i$. The $\mathrm{Aut}(E)$-module $\mathrm{Ext}^1(E,E)$ is then exactly the $\gl(\alpha)$-module $R(Q,\alpha)$ consisting of all representations of $Q$ of dimension $\alpha$ (we refer to (\textit{loc. cit.}) for the notations). This point of view identifies the quotient $\mathrm{Ext}^1(E,E)_0 \git \mathrm{Aut}(E)$ we have in mind with a closed subscheme of $R(Q,\alpha)\git \gl(\alpha)$, and (\textit{loc. cit.}) shows that the coordinate ring of the latter is generated by traces along oriented cycles in the quiver $Q$. But we also need a precise description of the relations between these generators (the \textit{second main theorem for invariant theory}). Once we have a convenient enough statement about these relations we can describe the completed local ring of $\mathcal{SU}_C(r)$ at $E$. 

When $r=3$ the decomposition (\ref{polystable}) ensures that there are only five cases to deal with, according to the values of the $r_i$'s and $\rho_i$'s.

\ssection The case of a stable bundle is obvious, and the case $r_1=2, r_2=1$ is a special case of the situation studied in \cite[III]{localstr}: $\mathcal{SU}_C(3)$ is \'etale isomorphic at $E$ to a rank $4$ quadric in $\mathbb A^9$. Here quivers do not provide a shorter proof. \mylabel{theorem}{easycase}

\ssection Let us look at the three other cases, where every $E_i$ in (\ref{polystable}) is invertible. The generic case consists of bundles $E$ which are direct sum of $3$ distinct line bundles. It has already been performed in \cite[V]{localstr}, but may also be recovered in a more convenient fashion as an easy consequence of \cite{LBP}: the generators of \cite[Lemma V.1]{localstr} then arise nicely as traces along closed cycles in the quiver \mylabel{theorem}{generic}

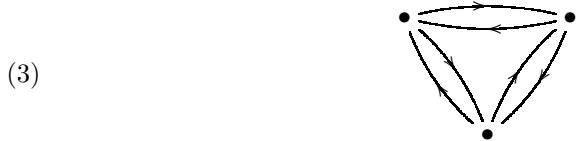
\begin{equation}\label{quiver}
\begin{array}{c}
\xymatrix@R=34pt@C=20pt{
\bullet \ar@{-}|@{>}@/^1.5mm/[rr] \ar@{-}|@{>}@/^1.5mm/[rd]& & \bullet \ar@{-}|@{>}@/^1.5mm/[ll] \ar@{-}|@{>}@/^1.5mm/[ld]\\
& \bullet \ar@{-}|@{>}@/^1.5mm/[ul] \ar@{-}|@{>}@/^1.5mm/[ur]&
}
\end{array}
\end{equation}

\noindent (note that there should be two loops on each vertex; but $\alpha=(1,1,1)$ implies that we can restrict ourselves to the quiver (\ref{quiver})). It is easy too to infer from (\ref{quiver}) the relation found by Laszlo; but, although \cite{LBP} gives a way to produce all the relations, this description turns out to be quite inefficient even in the present case (note however that, in order to conclude here, it is enough to remind that we know a priori the dimension of $\mathrm{Ext}^1(E,E) \git \mathrm{Aut}(E)$).

In the remaining two cases we already know that the tangent cone at $E$ must be a quadric (in $\mathbb A^9$) of rank $\leqslant 2$ (see \cite[V]{localstr}). We give now more precise statements.

\ssection\mylabel{theorem}{2by2} Suppose that $\rho_1=2$, i.e. that $E=(L\otimes V) \oplus L^{-2}$ where $L$ is a line bundle of degree $0$ with $L^3 \not\simeq \O$ and $V$ a vector space of dimension $2$. We have to consider here the ring of invariant polynomials on the representation space $R(Q,(2,1))$ of the quiver $Q$

\

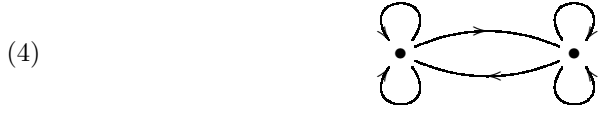
\begin{equation}\label{quiver2}
\begin{array}{c}
$$\xymatrix{
 \bullet \ar@{-}|@{>}@/^3mm/[rr]  \ar@(ur,ul)[] \ar@(dr,dl)[] &  & \bullet \ar@{-}|@{>}@/^3mm/[ll] \ar@(ul,ur)[] \ar@(dl,dr)[] 
}$$
\end{array}
\end{equation}

\

\

\noindent under the action of $\gl(V)\times\mathbb G_m$. Since the second vertex corresponds to a $1$-dimensional vector space it is enough to consider the quiver obtained by deleting the two loops on the right, and in fact we are brought to the action of $\gl(V)$ on $\mathrm{End}(V) \oplus \mathrm{End}(V) \oplus \mathrm{End}(V)^{\leqslant 1} \subset \mathrm{End}(V)^{\oplus 3}$, where $\mathrm{End}(V)^{\leqslant 1}$ denotes the space of endomorphisms of $V$ of rank at most $1$: this simply means that 
$$k[R(Q,(2,1))]^{\gl(V)\times\mathbb G_m}\simeq \left(k[R(Q,(2,1))]^{\mathbb G_m}\right)^{\gl(V)},$$
\noindent and that $k[R(Q,(2,1))]^{\mathbb G_m}$ gets naturally identified (as a $\gl(V)$-module) with $\mathrm{End}(V) \oplus \mathrm{End}(V) \oplus \mathrm{End}(V)^{\leqslant 1} \oplus k \oplus k$, the last two summands being fixed under the induced operation of $\gl(V)$. 

Let us now translate this discussion in a more geometric setting. Since (\ref{Ext}) identifies here $\mathrm{Ext}^1(E,E)$ with 
$$\left(H^1(C,\O)\otimes\mathrm{End}(V)\right)\oplus \left(H^1(C,L^{-3})\otimes V^\ast\right)\oplus\left(H^1(C,L^3)\otimes V\right)\oplus \left(H^1(C,\O)\otimes k\right),$$
\noindent we can identify the $\mathrm{Aut}(E)$-module $\mathrm{Ext}^1(E,E)_0$ with the $\gl(V)\times\mathbb G_m$-module 
$$\left(H^1(C,\O)\otimes\mathrm{End}(V)\right)\oplus \left(H^1(C,L^{-3})\otimes V^\ast\right)\oplus\left(H^1(C,L^3)\otimes V\right),$$ 
so that, up to the choices of some basis of the different cohomology spaces, any element of $\mathrm{Ext}^1(E,E)_0$ can be written $(a_1,a_2,\lambda,v) \in \mathrm{End}(V) \oplus \mathrm{End}(V) \oplus V^\ast \oplus V$. The map $(a_1,a_2,\lambda,v) \mapsto (a_1,a_2,a_3=\lambda \otimes v) \in \mathrm{End}(V)^{\oplus 3}$ identifies the quotient $\mathrm{Ext}^1(E,E)_0 \git \mathrm{Aut}(E)$ with the closed subscheme of $\mathrm{End}(V)^{\oplus 3} \git \gl(V)$ defined by the equation $\det a_3=0$. A presentation of the invariant algebra $k[\mathrm{End}(V)^{\oplus 3}]^{\gl(V)}$ can be found in \cite{drensky} (note that another presentation of this ring had been previously given in \cite{formanek}): if we let $b_i$ denote the traceless endomorphism $a_i-\frac{1}{2}\tr(a_i)\mathrm{id}$, this invariant ring is generated by the following ten functions
\begin{equation} \label{gen} 
\begin{array}{c}
 \displaystyle{u_i=\tr(a_i)\ \text{with $1 \leqslant i \leqslant 3$},\ v_{ij}=\tr(b_i b_j)\ \text{with $1\leqslant i \leqslant j \leqslant 3$},} \\
\\
 \displaystyle{w=\sum_{\sigma\in\mathfrak{S}_3} \varepsilon(\sigma) \tr(b_{\sigma(1)}b_{\sigma(2)}b_{\sigma(3)}),}
\end{array}
\end{equation}
\noindent subject to the single relation $w^2+18 \det(v_{ij})=0$. We have thus obtained the following result:
\begin{e-proposition}
If $E=\left(L \otimes V\right) \oplus L^{-2}$ with $L^3 \not\simeq \O$, then $\mathcal{SU}_C(3)$ is \'etale locally isomorphic at $E$ with the subscheme of $\mathbb A^{10}$ defined by the two equations 
$$X_{10}^2+18(X_4X_5X_6+2X_7X_8X_9-X_6X_7^2-X_5X_8^2-X_4X_9^2)=0$$
$$\text{ and } X_3^2-2X_6=0$$
\noindent at the origin. Its tangent cone is a double hyperplane in $\mathbb A^9$.
\end{e-proposition}
\ssection\mylabel{theorem}{3by3} Suppose now that $\rho_1=3$, i.e. that $E = L \otimes V$ where $V$ is a vector space of dimension $3$ (and $L$ a line bundle of order $3$). By the same argument as in \cite[Proposition V.4]{localstr} we know that the tangent cone at such a point is a rank $1$ quadric. But an explicit description of an \'etale neighbourhood is available, thanks to \cite{two3by3}. The space $\mathrm{Ext}^1(E,E)_0$ is isomorphic to $H^1(C,\O)\otimes \mathrm{End}_0(V)$ and, if we fix a basis of $H^1(C,\O)$, any of its element can be written $(x,y)\in \mathrm{End}_0(V)\oplus\mathrm{End}_0(V)$. The ring of invariants $k[H^1(C,\O) \otimes \mathrm{End}_0(V)]^{\gl(V)}$ is then generated by the nine functions $\tr(x^2)$,$\tr(xy)$, $\tr(y^2)$, $\tr(x^3)$, $\tr(x^2y)$, $\tr(xy^2)$, $\tr(y^3)$, $v=\tr(x^2y^2)-\tr(xyxy)$ and $w=\tr(x^2y^2xy)-\tr(y^2x^2yx)$;
moreover the ideal of relations is principal, generated by an explicit equation (see (\textit{loc. cit.})). 

As a result of this case-by-case analysis we conclude that $\mathcal{SU}_C(3)$ is a local complete intersection, as announced in the introduction.

\section{On the local structure of $\mathcal S$}

We know from \cite{angela} that the involution $\sigma$ associated to the double covering given by the theta map
$$\theta \colon \mathcal{SU}_C(3) \to |3\Theta|$$
\noindent acts by $E \mapsto \iota^\ast E^\ast$, where $\iota$ stands for the hyperelliptic involution. The local study of its ramification locus thus reduces to an explicit analysis of the behaviour of $\sigma$ through the \'etale morphisms resulting from Luna's theorem. Once again it comes  to a case-by-case investigation.

\ssection When $E$ is stable there is nothing to say. If $E=F \oplus L$ (with $F$ a stable bundle of rank $2$ and $L=(\det F)^{-1}$) we have to understand the action of the linearization of $\sigma$ on 
$$\mathrm{Ext}^1(E,E)_0\simeq \mathrm{Ext}^1(F,F) \oplus \mathrm{Ext}^1(F,L) \oplus \mathrm{Ext}^1(L,F)$$
\noindent (note that we tacitly identify $\mathrm{Ext}^1(F,F)$ with its image in $\mathrm{Ext}^1(F,F) \oplus H^1(C,\O) \subset \mathrm{Ext}^1(E,E)$ by the map $\omega \mapsto (\omega,-\tr (\omega))$).

Since $\sigma(E)=E$, $\iota^\ast F^\ast$ must be isomorphic to $F$, and $\sigma$ identifies $\mathrm{Ext}^1(F,L)$ and $\mathrm{Ext}^1(L,F)$; let us choose a basis $X_1,X_2$ of $\mathrm{Ext}^1(F,L)^\ast$, and call $Y_1,Y_2$ the corresponding basis of $\mathrm{Ext}^1(L,F)$. We need here to recall precisely from \cite{localstr} the explicit description of the coordinate ring of $\mathrm{Ext}^1(E,E)_0\git \mathrm{Aut}(E)$ mentionned in \ref{easycase}: it is generated by $k[\mathrm{Ext}^1(F,F)]$ and the four functions $u_{ij}=X_iY_j$, subject to the relation $u_{11}u_{22}-u_{12}u_{21}=0$. 

It follows from our choice that $\sigma$ maps $u_{ij}$ to $u_{ji}$. Furthermore we claim that $\sigma$ acts identically on $\mathrm{Ext}^1(F,F)$: as a stable bundle, $F$ corresponds to a point of the moduli space $\mathcal{U}(2,0)$, whose tangent space is precisely isomorphic to $\mathrm{Ext}^1(F,F)$. The action of $\sigma$ on this vector space is the linearization of the one of $F \in \mathcal U(2,0) \mapsto \iota^\ast F^\ast$. Using that $\mathcal U(2,0)$ is a Galois quotient of $J_C \times \mathcal{SU}_C(2)$, our claim comes from the fact that $\sigma$ is trivial on both $J_C$ and $\mathcal{SU}_C(2)$.

Since the coordinate ring of the fixed locus of $\sigma$ in $\mathrm{Ext}^1(E,E)_0\git \mathrm{Aut}(E)$ is the quotient of the one of $\mathrm{Ext}^1(E,E)_0\git \mathrm{Aut}(E)$ by the involution induced by $\sigma$ we may conclude that $\mathcal S$ is \'etale locally isomorphic at $E$ to the quadric cone in $\mathbb A^8$ defined by $X_3^2-X_1 X_2=0$.

\ssection Consider now the situation of \ref{generic}: let us write $E=L_1\oplus L_2 \oplus L_3$ with $L_i \not\simeq L_j$ if $i\neq j$. We have $\mathrm{Ext}^1(E,E) \simeq \bigoplus_{i,j} \mathrm{Ext}^1(L_i,L_j)$; let us choose for $i \neq j$ a non-zero element $X_{ij}$ of $\mathrm{Ext}^1(L_i,L_j)^\ast$ such that $X_{ji}$ corresponds to $X_{ij}$ through the isomorphism $\mathrm{Ext}^1(L_i,L_j) \simeq \mathrm{Ext}^1(L_j,L_i)$ induced by $\sigma$ and the natural isomorphisms $\iota^\ast L_i^\ast \simeq L_i$. It then follows from \ref{generic} (see \cite{localstr} for a complete proof) that the ring $k[\mathrm{Ext}^1(E,E)_0]^{\mathrm{Aut}(E)}$ is generated by $k[\ker(\bigoplus_i \mathrm{Ext}^1(L_i,L_i) \to H^1(C,\O))]$ and the five functions
$Y_1=X_{23}X_{32}$, $Y_2=X_{13}X_{31}$, $Y_3=X_{12}X_{21}$, $Y_4=X_{12}X_{23}X_{31}$, $Y_5=X_{13}X_{32}X_{21}$,
subject to the relation $Y_4Y_5-Y_1Y_2Y_3=0$. One easily checks that the involution $\sigma$ fixes $k[\ker(\bigoplus_i \mathrm{Ext}^1(L_i,L_i) \to H^1(C,\O))]$, $Y_1$, $Y_2$ and $Y_3$, while it sends $Y_4$ to $Y_5$. The fixed locus $\mathrm{Fix}(\sigma)$ is then defined by the equation $Y_4-Y_5=0$, so that $\mathcal S$ is \'etale locally isomorphic to the hypersurface in $\mathbb A^8$ defined by $Z_4^2-Z_1 Z_2 Z_3=0$.
Its tangent cone is a double hyperplane.

\ssection In the situation of \ref{2by2} we have to make a more precise choice of the non-zero elements of $\mathrm{Ext}^1(L^{-2},L)$ and $\mathrm{Ext}^1(L,L^{-2})$, so as to make them correspond through $\sigma$ and the natural isomorphism $\iota^\ast L^\ast \simeq L$; such a choice ensures that $\sigma$ operates on $\mathrm{Ext}^1(E,E)_0$ in the following way:
$$(x,y,\lambda,v)\in \mathrm{End}(V)^{\oplus 2} \oplus V{}^\ast \oplus V\mapsto ({}^tx,{}^ty,{}^t v,{}^t\lambda),$$
\noindent so that we know how it acts on the generators of $k[\mathrm{Ext}^1(E,E)_0 \git \mathrm{Aut}(E)]$ given in (\ref{gen}): $\sigma$ fixes $u_i$, $v_{ij}$, and sends $w$ to $-w$. This implies that the fixed locus is defined by the equation $w=0$. The sextic $\mathcal S$ is \'etale locally isomorphic to the subscheme of $\mathbb A^9$ whose ideal is generated by the two equations 
$$X_4X_5X_6+2X_7X_8X_9-X_6X_7^2-X_5X_8^2-X_4X_9^2=0\ \text{and}\  X_3^2-2X_6=0;$$
\noindent its tangent cone is therefore the cubic hypersurface of $\mathbb A^8$ defined by $2 X_7 X_8 X_9 -X_5 X_8^2-X_4 X_9^2=0$.

\ssection We are now left with the last case, where $E$ is of the form $L \otimes V$ (with $L^3=\O$): $\mathrm{Ext}^1(E,E)_0$ is then isomorphic to $H^1(C,\O) \otimes \mathrm{End}_0(V)$, and $\sigma$ acts by $\omega\otimes a \in  H^1(C,\O) \otimes \mathrm{End}_0(V) \longmapsto \omega \otimes {}^t a$.
This induces an action on $k[H^1(C,\O)\otimes\mathrm{End}_0(V)]^{\mathrm{Aut}(E)}$ which fixes the first eight generators of \ref{3by3}, and acts by $-1$ on the last one, namely $w$; the fixed locus is thus defined in $\mathrm{Ext}^1(E,E)_0 \git \mathrm{Aut}(E)$ by the linear equation $w=0$. 

The sextic $\mathcal S$ is then \'etale locally isomorphic to an hypersurface in $\mathbb A^8$ defined by an explicit equation; writing down this equation shows that its tangent cone is a triple hyperplane.


\begin{thebibliography}{00}
% please try to use the bibitem system -
% the references should be in alphabetical order of authors' names.
% Articles with a single author first, author will 1 co-author next,
% then author with several co-authors;


% \bibitem{label}
% Text of bibliographic item

\bibitem{two3by3} H. Aslaksen, V. Drensky, L. Sadikova, \textit{Defining relations of invariants of two {$3\times3$} matrices}, J. Algebra \textbf{298} (2006) 41--57.

\bibitem{drensky} V. Drensky, \textit{Defining relations for the algebra of invariants of {$2\times
              2$} matrices}, Algebr. Represent. Theory \textbf{6} (2003) 193--214.

\bibitem{formanek} E. Formanek, \textit{Invariants and the ring of generic matrices}, J. Algebra \textbf{89} (1984) 178--223.
 
\bibitem{localstr} Y. Laszlo, \textit{Local structure of the moduli space of vector bundles over curves}, Comment. Math. Helv. \textbf{71} (1996) 373--401.

\bibitem{LBP} L. Le Bruyn, C. Procesi, \textit{Semisimple representations of quivers}, Trans. Amer. Math. Soc. \textbf{317} (1990) 585--598.

\bibitem{minh} Q. M. Nguy$\tilde{\hat{\text{e}}}$n, \textit{Vector bundles, dualities, and classical geometry on a curve of genus two}, Preprint \texttt{arXiv:math.AG/0702724}, to appear in Internat. J. Math.

\bibitem{angela} A. Ortega, \textit{On the moduli space of rank 3 vector bundles on a genus 2 curve and the {C}oble cubic}, J. Algebraic Geom. \textbf{14} (2005) 327--356.

\bibitem{orth} O. Serman, \textit{Moduli spaces of orthogonal bundles over an algebraic curve}, Preprint \texttt{arXiv:math.AG/0609520}.

\end{thebibliography}
\end{document}